\documentstyle[12pt,euscript,leqno]{amsart}
\newtheorem{thm}{Theorem}[section]
\newtheorem{lem}[thm]{Lemma}
\newtheorem{cor}[thm]{Corollary}

\numberwithin{equation}{section}

\newcommand{\Om}{\Omega}

\newcommand{\RR}{{\mathbb{R}}}

\renewcommand{\phi}{\varphi}

\newcommand{\UU}{{\mathcal U}}

\begin{document}

\title[ON THE LOCAL THEORY OF PRESCRIBED JACOBIAN EQUATIONS]
{ON THE LOCAL THEORY OF PRESCRIBED JACOBIAN EQUATIONS}
\subjclass[2010]{Primary 35J66,35J25; Secondary  78A05. }
\author{Neil S. Trudinger${}^\dagger$}
\thanks{${}^\dagger$Research supported by Australian
Research Council Grant.}
\address{${}^\dagger$Centre for 
Mathematics and its Applications,
Australian National University, Canberra, ACT 0200, Australia.}
\email{neil.trudinger@@anu.edu.au}

\maketitle

\vspace{0.3cm}

\noindent { Abstract.}
\vspace{0.3cm}

We develop the fundamentals of a local regularity theory for prescribed Jacobian equations which
extend the corresponding results for optimal transportation equations. In this theory the cost function is extended to a \emph{generating function} through dependence on an additional scalar variable.
In particular we recover in this generality the local regularity theory for potentials of Ma, Trudinger and Wang, along with the subsequent development of the underlying convexity theory.

\vspace{1cm}

%
%
\section{Introduction}
Let $\Om$ be a domain in Euclidean $n$-space, $\mathbb{R}^{\it n}$, and
$Y$ a mapping from $\Om \times\mathbb{R} \times\RR^{\it n}$ into $\mathbb{R}^{\it n}$. The
\emph{prescribed Jacobian equation}, PJE, is a partial differential equation of
the form,
%
%
\begin{equation} \label{def:PJE}
 \det DY (\cdot,u,Du)  =  \psi(\cdot,u,Du),
\end{equation}

\noindent where $\psi$ is a given scalar function on
$\Om\times\RR\times{\RR}^n$ and $Du$ denotes the gradient of the
function $u:\Om \to\RR$. Denoting points in  $\Om \times\RR \times{\RR}^n$
by $(x,z,p)$, we see that the special case,
%
%
\begin{equation} \label{def:grad}
 Y(x,z,p) = p,
\end{equation}

\noindent corresponds to the standard Monge-Amp\`ere equation,
%
%
\begin{equation} \label{def:MA}
\det D^2u =  \psi(\cdot,u,Du).
\end{equation}

\noindent We will always assume that the matrix $Y_p$ is invertible,
that is $\det Y_p \neq 0$, whence we may write \eqref{def:PJE} as a general
equation of Monge-Amp\`ere type,
%
%
\begin{equation} \label{def:MAE}
\det [D^2u - A(\cdot,u,Du) ] =  B(\cdot,u,Du),
\end{equation}

\noindent where
%
%
\begin{equation} \label{def:AB}
A = - Y^{-1}_p(Y_x+Y_z\otimes p),
\quad B = (\det Y_p)^{-1}\psi.
\end{equation}

\noindent A function $u\in C^2(\Om)$ is degenerate elliptic,
(elliptic), for equation \eqref{def:MAE}, henceforth called \emph{admissible},
whenever
%
%
\begin{equation} \label{def:elliptic}
D^2u - A(\cdot,u,Du) \ge 0, \quad (>0),
\end{equation}
\noindent in $\Om$. If $u$ is an elliptic solution of \eqref{def:MAE}, then the function $B(\cdot,u,Du)$
is positive. Accordingly we will assume throughout that $B$ is at least non-negative
in $\Om\times\RR\times\RR^n$, that is $\psi$ and $\det Y_p$ have the same
sign.

 The \emph{second boundary
value problem} for the prescribed Jacobian equation is  to prescribe the image,
%
%
\begin{equation} \label{def:bvp}
Tu(\Om): = Y (\cdot,u,Du)(\Om) = \Om^*,
\end{equation}

\noindent where $\Om^*$ is another given domain in
 $\RR^n$. When
$\psi$ is separable, in the sense that
%
%
\begin{equation} \label{def:separable}
|\psi(x,z,p)| = f(x)/g\circ Y(x,z,p),
\end{equation}
\noindent for positive $f, g\in L^1(\Om)$, $L^1(\Om^*)$ respectively, then a
necessary condition for the existence of an admissible solution,
for which the mapping $T$ is a diffeomorphism, to
the second boundary value problem \eqref{def:PJE}, \eqref{def:bvp} is the \emph{mass
balance} condition,
%
%
\begin{equation} \label{def:massbalance}
\int_\Om f=\int_{\Om^*} g.
\end{equation}
\vspace{0.3cm}

\noindent For the standard  Monge-Amp\`ere equation \eqref{def:MA}, with $Tu = Du$,
 the classical solvability of the second boundary value problem, under
 the mass balance condition \eqref{def:massbalance}, was proved by Delan\"oe, ($n=2$), \cite{D},
 Caffarelli \cite{C2} and Urbas\cite{U}, under the hypothesis that both domains,
 $\Om$ and $\Om^*$ are uniformly convex. As already pointed out in
 \cite{LTU}, \eqref{def:bvp} implies a nonlinear \emph{oblique} boundary condition.
A weaker interpretation of  the boundary condition \eqref{def:bvp} arises through optimal
transportation, in which case Caffarelli \cite{C1} proved that the convexity of
the target $\Om^*$ suffices for local smoothness of solutions.

Interest in the general case was  stimulated in the last decade through its application
 to  regularity in optimal transportation, (\cite{MTW,TW2}).  Here we are given
  a cost function $c:\RR^n \times \RR^n \to \RR$ and the vector field $Y$
  is generated by the equation,
  %
  %
\begin{equation} \label{def:costfunction}
c_x(x,Y(x, p))=  p,
\end{equation}
which we assume to be uniquely solvable for $p\in c_x(\Om\times\Om^*)$,
with non-vanishing determinant, that is $\det  c_{x, y}(x, y)\ne 0$, for
all $(x, y)\in \Om\times \Om^*$. In the corresponding Monge-Ampere
equation \eqref{def:MAE}, we have
%
%
\begin{equation} \label{def:OTE}
A(x, z, p)  = c_{xx}(x,Y(x, p)),\quad B(x, z, p) =\det c_{x,y}(x,Y (x,p)\psi (x, z, p).
\end{equation}

\noindent Conditions for local regularity were found in \cite{MTW} and
for global regularity
in \cite{TW3}, with the latter being extended to general prescribed Jacobian
equations in \cite{T1,T2}. 

In this article we consider a more general situation where the cost function is
replaced by a smooth \emph{generating function} $G:\RR^n \times \RR^n\times\RR \to \RR$,
with the resultant vector field $Y$ also depending on $u$.
This enables the corresponding local theory to embrace
recent work in near field optics, \cite{KO,KW}. The assumptions on the generating
 function $G$ are parallel
to those introduced for cost functions in optimal transportation in \cite{MTW,TW3}.
Accordingly we let
$\UU$ be an open set in  $ \RR^n\times\RR^n$ and $I$ an open interval in
$\RR$. For points
$(x, y)\in\UU$, we denote their corresponding projections by
$\UU^*_x=\{y\in \RR^n\ |\ (x, y)\in \UU\}$, $\UU_y=\{x\in \RR^n\ |\ (x,
y)\in \UU\}$ and write $\UU^{(1)}=\bigcup\{\UU_y\ |\ y\in\RR^n\}$ and
$\UU^{(2)}=\bigcup\{\UU^*_x\ |\ x\in\RR^n\}$ Denoting points in $I$ by $z$, 
we assume that
$G$ is smooth in $\UU\times I,  G_z \ne 0$ and
\vspace{0.2cm}
\begin{itemize}
\item[{\bf G1}:]
For each $(x,y) \in \UU$, there exists an open interval
$I(x,y)\subset I$ such that the mapping $(G_x,G)(x,\cdot,\cdot) $ is one-to -one
in $y\in\UU^*_x, z\in I(x,y)$, for each $x\in\UU^{(1)}$.
\vspace{0.2cm}

\item[{\bf G2}:]
For each $(x,y)\in\UU, z\in I(x,y)$,
$\det E(x,y,z)\neq 0$, where $E$ is the $n\times n$ matrix given by

%
%
\begin{equation} \label{def:E}
E = [E_{x,y}] = [ G_{x,y} - (G_z)^{-1}G_{x,z}\otimes G_y].
\end{equation}

\end{itemize}
\vspace{0.3cm}
From G1 and G2, the vector field $Y$, together with a scalar function $Z$,
are generated by G through the equations,

%
%
\begin{equation} \label{def:YZ}
G_x( x, Y, Z) = p,\quad G(x,Y,Z) =u.
\end{equation}
\vspace{0.3cm}

\noindent The significance of the additional function $Z$ will become apparent later.
Note that the Jacobian determinant of the mapping $(y,z) \to  (G_x,G) (x,y,z)$
is  $G_z \det E, \neq 0$ by G2, so that $Y $and $Z$ are accordingly smooth.
Also by differentiating \eqref {def:YZ}, with respect to $p$, we obtain $Y_p = E^{-1}$.
Next using \eqref{def:AB} or differentiating \eqref{def:YZ} for $p = Du$, with respect to $x$, 
we obtain that
the prescribed Jacobian equation \eqref{def:PJE}  is a Monge-Amp\`ere equation of the
form \eqref {def:MAE}  with

%
%

\begin{align}\label{def:GPJE}
A(x, u, p)  &= G_{xx}[x,Y(x,u, p),Z(x,u,p)] , \\
 B(x, u, p) &=\det E(x,Y,Z)\psi (x, u, p) \notag
\end{align}

\noindent and is well defined in domains $\Om\in\UU^{(1)}$ for $Y\in\UU^*_x, Z\in I(x,Y), x\in\Om$.
Note that the latter restrictions may automatically place constraints on $u$ and $Du$.

In the optimal transportation case, we have
%
%

\begin{equation}\label{def:OTG}
G(x,y,z) = c(x,y) - z,\quad G_z = -1,\quad I = I(x,y) = \RR,
\end{equation}
\noindent and we recover \eqref{def:costfunction}  and 
\eqref{def:OTE}  above.

Note also that by adjusting the dependence of $G$ on $z$, we can fix the interval $I$, as well
as the sign of $G_z$, as we wish.
For convenience with other examples, we will assume either $I =(0,\infty)$ or $I = \RR$ as above
and assume $G_z < 0$. Let us also denote 
$$ \Gamma = \Gamma(\UU) = \big\{(x,y,z)\in\UU\times I\mid z\in I(x,y)\big\}$$

We mention also that the idea behind conditions G1 and G2 is to determine the mapping $T$
from a tangential intersection of the graphs of the functions $u$ and $G (\cdot,y,z)$. If the graph
of $G$ is also a local support from below, we obtain the ellipticity condition (1.6), which
corresponds to a local convexity. The geometric picture is further amplified in \cite{T3},
where the theory developed here is already outlined. Knowledge of the defining function $G$
can also lead to a more efficient derivation of the  Monge-Amp\`ere equation, using (1.14)
rather than computing it directly from (1.5). Also note that when the graph of $G$ is a local support from above we obtain the complementary ellipticity condition,

%
%

\begin{equation}\label {def:c-e}
D^2u - A(\cdot,u,Du) \le 0,
\end{equation}

\noindent corresponding  to a local concavity. In this case by setting $u^- = -u, A^-(\cdot,u,p) = -A(\cdot,-u,-p)$, we also obtain a degenerate elliptic equation of the form (1.4) for the function $u^-$. Furthermore if $n=2$ and  $B > 0$ in (1.4), then either of the strict versions of (1.6) or (1.16) must hold.

For the regularity of weak solutions we need a global convexity theory. To develop this, 
we in turn need a dual condition,
which in the optimal transportation case is simply obtained by interchanging $x$ and $y$.
Using the property that the generating function $G$ is strictly monotone with respect
to $z$, we introduce a dual function $G^* = H$  on $\UU\times \RR$ by

%
%

\begin{equation}\label {def:H}
G[x,y,H(x,y,u)] = u.
\end{equation}

Clearly $H$ is well defined whenever $(x,y)\in\UU$ and $u\in G(x,y,\cdot)(I)$.
Furthermore we have the relations
%
%

\begin{equation}\label{def:DH}
H_x  = - G_x/G_z,  H_y = -G_y/G_z,  H_u = 1/G_z.
\end{equation}
This motivates the following dual condition:

\vspace{0.2cm}
\begin{itemize}
\item[{\bf G1*}:]

The mapping $Q: = -G_y/G_z$ is one-to-one in $x$, for all $y\in\UU^{(2)}, z \in I(x,y).$

\end{itemize}
\vspace{0.2cm}
\noindent Writing $J(x,y) = G (x,y,\cdot)I(x,y)$, we see that condition G1* corresponds to
G1 with $x$ and $y$ interchanged and $I(x,y)$ replaced by $J(x,y)$. Furthermore
the Jacobian matrix of the mapping $x\to Q(x,y,z)$ is $-E/G_z$ so its determinant is
automatically non-zero when condition G2 holds. Analogously to $\Gamma$ above, we
may also denote the dual set,

$$\Gamma^* = \big\{(x,y,u) \in\UU\times\RR\mid u\in J(x,y)\big\}.$$

Our next conditions extend the conditions A3 and A3w introduced for regularity in \cite{MTW,TW2}
and are expressed in terms of the matrix function $A$ in \eqref{def:MAE}, namely:

\vspace {0.3cm}

\noindent {\bf G3}\ ({\bf G3w}) \vspace{-0.5cm}
$$A^{kl}_{ij}\xi_i\xi_j\eta_k\eta_l: = (D_{p_kp_l}A_{ij}) \xi_i\xi_j\eta_k\eta_l > (\ge)\ 0,
        $$

\vspace {0.2cm}

 for all $(x,Y)\in \UU, Z\in I(x,Y),  \xi,\eta \in \RR^n$ such that $\xi.\eta = 0$.

\vspace {0.3cm}

\noindent As in \cite{T1}, we may write equivalently that $A$ is \emph{strictly regular}, 
(\emph{regular}), in the set 

$$\Gamma^\prime = \big \{(x,u,p) \mid x\in\UU^{(1)}, u = G (x,y,z), p = G_x(x,y,z), 
\text{for some}\ y\in\UU^*_x, z\in I(x,y)\}.$$
\vspace{0.3cm}
For our convexity theory we will also need  a condition
that the matrix function $A$ is monotone with respect to
$u$, that is

\vspace{0.3cm}
\noindent {\bf G4}\ ({\bf G4w}) \vspace{-0.5cm}
$$D_uA_{ij}\xi_i\xi_j > (\ge)\ 0,$$

\vspace {0.2cm}

  for all $ (x,Y)\in \UU, Z\in I(x,Y),  \xi \in \RR^n$.

\vspace {0.3cm} 

\noindent In this paper we will not use the strict monotonicity G4. 

\vspace{0.3cm}
We illustrate the above conditions with the example of a parallel beam from \cite{LT2}.
For
 $\UU = \RR^n\times\RR^n$ and $I = (0,\infty)$, we define
%
%
\begin{equation}\label  {def:example}
G(x,y,z) = {1\over 2z} - {z\over 2}|x-y|^2.
\end{equation}

\noindent Then $G$ satisfies G1,G2,G1*,G3, G4 with
\vspace{0.3cm}
\begin{align}
G_z &= -{1\over 2}(z^{-2}  +|x-y|^2) < 0,\quad G_x = z(x-y),\quad G_y = z(y-x),\notag \\
E& = z\big\{(1 + z^2|x-y|^2)I - z^2(x-y)\otimes(x-y)\big\}/(1 + z^2|x-y|^2)),\notag\\
\quad \text{det}E& = z^n(1 - z^2|x-y|^2)/(1+z^2|x-y|^2) > 0, z\in I(x,y),\notag\\
Y& = x +{2uDu\over(1- |Du|^2)}, \quad Z = {(1- |Du|^2)\over 2u}, \quad A = {-ZI},\notag\\
H& = {1\over u + (u^2 + |x-y|^2)^{1/2}} ,\quad I(x,y) = (0,{1\over|x-y|}), J(x,y) = (0,\infty).\notag
\end{align}

\vspace{0.3cm}
\noindent The corresponding Monge-Amp\`ere equation,
%
%
\begin{equation}
\text{det}\big\{ D^2u + {(1 - |Du|^2)\over 2u}I\big\} = {(1 - |Du|^2)^{n+1}\over(1+|Du|^2)(2u)^n}\psi,
\end{equation}
\noindent
is well defined for $u > 0$ and $|Du| < 1$.

In the underlying physical model, a parallel beam of light directed upwards through $\RR^n$
is reflected, in accordance with Snell's law, from the graph of u to a target back in $\RR^n$.
The restrictions $u > 0$ and $|Du| < 1$ are thus obvious and  $Tu(x)$ is the point in  $\RR^n$
reached by a incoming ray in the upwards direction  through the point $x$. The graph of $G$ in $\RR^{n+1}$ is a focal paraboloid in that every vertical ray is reflected to the focus $(y,0)$. When
$\psi = f/g(T)$, (1.20) is the PDE satisfied by a function $u$ for which the reflection map $Tu$
pushes forward the density $f$ to the density $g$.
\vspace{0.3cm}

In the next section, we show that a convexity theory, with respect to generating functions,
replicating usual convexity, and more
generally the optimal transportation case as developed in \cite{KM,Loe,TW2,TW4,V2} can be built under conditions G1, G2, G1*,G3w and G4w.
In particular we prove that the local convexity of smooth functions implies their global 
convexity for appropriately convex domains, that 
normal mappings are determined by
sub-differentials and that sections and contact sets are convex in the extended sense.
 As well we prove that 
domain convexity with respect to generating functions is a special case of that determined
by vector fields in \cite{T1,T2}. Our  proofs follow the optimal transportation case, as presented
for example in \cite{TW2} ,
although they also depend on more intricate calculations.

In Section 3, we derive  the dual Monge-Amp\`ere equation, 
satisfied by the $G$-transform $v$ of a $G$-convex solution $u$ of (1.1), given by
 %
%
\begin{equation} \label{def:funcv}
v(y) = u^*_G(y) = \text{sup}_{\Om} H(\cdot,y,u).
\end{equation}
and prove that conditions G3 and G3w are invariant under duality.

Section 4 is devoted to the existence and regularity of generalized solutions of the second
boundary value problem (1.7), (1.8), with initial and target domains
satisfying $\Om\times\Om^*\subset\UU$. For existence we follow the approach in \cite{CO,KO}, which
corresponds to the existence of potential functions  in optimal transportation.
Here we  need an additional condition to control gradients of solutions,
which we may express as:
 
\begin{itemize}
\item[{\bf G5}:]
There exists  constants $m_0\ge-\infty, K_0 \ge 0$, such that $(m_0,\infty) \subset J(x,y)$ 
and

$$  |G_x(x,y,z)| \le K_0 $$

\noindent for all $x\in\Om, y\in\Om^*, G(x,y,z) > m_0$
\end{itemize}
\vspace{0.3cm}
The example (1.18) clearly satisfies G5 for  $m_0 =0, K_0 = 1$, while in the optimal transportation
case (1.15), $m_0 = -\infty$. Note that by translation, we may always assume either $m_0 = 0$
or $m_0 = -\infty$. Under condition G5, there exists a generalized solution $u$
with $u(x_0) = u_0$, provided $u_0 > m_0 + K_0d(x_0,\partial\Om)$. More generally
 we can also replace $G$ by composites $\mu(G)$ for suitable
smooth functions $\mu \in C^1(m_0,\infty), \mu^\prime \ne 0$.  This enables us to embrace
more examples such as near field reflection from a point source, as treated in \cite{KO,KW};
(see Section 4).

Following \cite{MTW}, we then prove local regularity under conditions
G1,G2,G1*G3 and G4w for target domains satisfying the appropriate convexity conditions. 
We remark that the existence of globally smooth elliptic solutions under corresponding
conditions, including stronger domain convexity conditions, follows from
the theory of the general prescribed Jacobian equation in \cite{T2} and its adaptation
 to near field reflection problems in \cite{LT2}.

\vspace{1cm}

%
%

\section{Convexity theory}

We begin with the appropriate definitions of convexity with respect to generating functions.

Let $\Om$ be a bounded domain in $\UU^{(1)}$, denote $\UU^*_\Omega = \cap_{x\in\Omega}\UU^*_x$
 and let $G$ be a generating function on $\UU\times I$, 
satisfying conditions G1 and G2 , with $I$ an open interval in
$\RR$ as in the previous section. A function $u\in C^0(\Om)$ is called $G$-\emph{convex} in $\Om$,
 if for each $x_0\in\Om$, there exists $y_0\in\UU^*_\Omega$ and $z_0\in\ I(x_0,y_0)$
 such that
 %
%
\begin{align}\label{def:Gconvex}
u(x_0) &= G(x_0,y_0,z_0), \\
 u(x) &\ge G(x,y_0,z_0)\notag
  \end{align}
\noindent for all $x\in\Om$. If $u$ is differentiable at $x_0$, then $y_0 = Tu (x_0) = 
Y(x_0,u(x_0),Du(x_0))$,
while if $u$ is twice differentiable at $x_0$, then

$$ D^2u(x_0) \ge D^{2}_xG(x_0,y_0,z_0),$$
\noindent that is $u$ is admissible for equation (1.4) at $x_0$. If $u\in C^2(\Om)$,
we call u \emph{locally} $G$-\emph{convex} in $\Om$ if this inequality holds for all $x_0\in\Om$.
We will also refer to functions of the form $G(\cdot,y_0,z_0)$ as \emph{$G$-affine} and as
a \emph{$G$-support} at $x_o$ if \eqref{def:Gconvex} is satisfied.

As in the optimal transportation case \cite{MTW}, we also have corresponding notions of 
domain convexity. There are various possibilities depending on what quantities are fixed.
For our purposes here we make the following definitions for domains $\Om$ and $\Om^*$ satisfying $\Om\times\Om^* \subset \UU$. 
\vspace{0.3cm}

\noindent The domain $\Om$ is $G$-\emph{convex} with respect to $y_0\in\UU^*_\Omega,\, z_0\in I (\Om,y_0)  = \cap_\Om I(\cdot,y_0)$ if
the image $Q_0(\Om): = -G_y/G_z (\cdot,y_0,z_0)(\Om)$ is convex in $\RR^n$;  
\vspace{0.3cm}

\noindent The domain $\Om^*$ is $G^*$-\emph{convex}
with respect to $(x_0,u_0)$, where  $x_0\in\UU_\Om^*$ and  
$u_0 \in J(x_0,\Om^*) =  \cap_{\Om^*}J(x_0,\cdot)$, if
the image $P_0(\Om^*): = G_x [x_0,\cdot,H(x_0,\cdot,u_0)](\Om^*)$ is convex in $\RR^n$.
\vspace{0.3cm}

Note when we use condition G3w below for convexity results and their consequences, we will assume at least that
the convex hulls of the images $Q_0(\Om)$ and $P_0(Tu(\Om))$
lie in $Q (\Gamma), G_x(\Gamma)$ respectively.
\vspace{0.2cm}

The second definition is clearly a special case of the notion of $Y^*$-convexity in \cite{T1},
as it is equivalent to the set $P_0(\Om^*) = \{p\in\RR^n\mid Y(x_0,u_0,p)\in\Om^*\}$ being convex.
Moreover we will define the domain $\Om^*$ to be $G^*$-convex with respect
 to a function $u \in C^0(\Om)$
if $\Om^*$ is $G^*$-convex with respect to each point on the graph of $u$, that is the sets $P_x(\Om^*)
 =  \{p\in\RR^n\mid Y(x,u(x),p)\in\Om^*\}$ are convex for each $x\in \Om$ The
relationship of the first definition and the notion of $Y$-convexity in \cite{T1} is treated in Lemma 2.4 below.
Our main Lemmas 2.1, 2.2 and 2.3 extend corresponding results for optimal transportation in 
\cite{TW3,KM,Liu} and \cite{FKM}.

  \begin{lem} 
Assume G1,G2,G1*,G3w and G4w  hold in $\UU\times I$ and that $u\in C^2(\Om)$
is locally $G$-convex in $\Om$. Then if $\Om$ is $G$-convex with respect to each point in
$(Y,Z)(\cdot,u,Du)(\Om)$,  $u$ is $G$-convex in $\Om$.
\end{lem}

More generally, we have for any $x_0\in\Om, y_0, z_0 = Y, Z(x_0,u(x_0),Du(x_0))$, the
$G$-affine function $G(\cdot ,y_0,z_0)$ is a $G$-support, provided $\Om$ is $G$-convex
with respect to $(y_0,z_0)$.

We defer the proof of Lemma 2.1, until the end of the section and proceed now to consider the
corresponding extensions of normal mappings and sections.
\vspace{0.3cm}

Let $u\in C^0(\Om)$ be $G$-convex in $\Om$. We define the \emph{$G$-normal} mapping
of $u$ at $x_0\in\Om$ to be the set:

$$Tu(x_0) = \big\{ y_0\in\UU_\Omega \mid 
 u(x) \ge G\big(x,y_0,H(x_0,y_0,u_0)\big)\text{ for all } x\in\Om\big\},$$
 
\noindent where $u_0 = u(x_0)$. Clearly $Tu$ agrees with our previous terminology when $u$ is differentiable
  and moreover in general
 
 $$ Tu (x_0) \subset Y(x_0,u(x_0),\partial u(x_0)),$$
 
 \noindent where $\partial u$ denotes the subdifferential of $u$.
 
\begin {lem}
Assume G1,G2,G1*,G3w,G4w  hold in $\UU\times I$ and suppose
$u\in C^0(\Om)$ is $G$-convex in $\Om$. Then for
any $x_0\in\Om$, we have

 $$ Tu(x_0)  = Y(x_0,u(x_0),\partial u(x_0)).$$
 
 \end{lem}
 
 For  a fixed $y_0\in Tu(x_0)$, corresponding $z_0 = H(x_0,y_0,u_0)$ and
 $\sigma > 0$, we define the \emph{$G$-section}, $S_\sigma$ by
 
 $$ S_\sigma  =  S_\sigma(x_0,y_0) = \big\{ x\in\Om \mid u(x) < G(x,y_0,z_0) + \sigma) \big\}$$
 
 \noindent and lower contact set $S_0$ by
 
 $$ S_0 = S_0(x_0,y_0) = \big\{ x\in\Om \mid u(x) = G(x,y_0,z_0) \big\}.$$
 
 \begin{lem}
Assume G1,G2,G1*,G3w,G4w  hold in $\UU\times I$ and suppose
$u\in C^0(\Om)$ is $G$-convex in $\Om$, with $\Om$ itself being $G$-convex
with respect to $y_0\in Tu(x_0)$ and $z_0 = H(x_0,y_0,u(x_0))$. Then the sets
$S_\sigma$ and $S_0$ are also $G$-convex with respect to $(y_0,z_0)$.
\end{lem}

The proofs of Lemmas 2.1,2.2 and 2.3 reduce to calculations that are more complicated versions of those which underwrite
our starting point in \cite{TW3,TW4}, relating domain convexity to the Monge-Amp\`ere
equation (1.4) . We will adopt the following notation:
%
%
\begin{align} \label{def:Notation}
  G_{i,\cdots , j,\cdots , z,\cdots, z}  &= 
  {\partial \over \partial x_i}  \, \cdots  
    {\partial \over \partial y_j}  \, \cdots  {\partial \over \partial z}\, \cdots
     {\partial \over \partial z} \, G , \\
 E_{i,j} &= E_{x_i,y_j}  \notag \\
 \big [ E^{i,j} \big ] &= E^{-1} = \big [ D_{p_j}Y^i \big ] ,\notag
 \end{align}
 
 \noindent so that in particular we have the following formulae for differentiation
 with respect to the $p$ variables:
 %
 %
 \begin{align} \label{def:D_p}
  Z_p &= - Y_p \; {G_y \over G_z} \\
  D_p &= Y_p \big (D_y- \;{G_y \over G_z} \, D_z \big ) \notag\\
  D_{p_k} &= E^{r,k}  \big ( D_{y_r} - \, {G_{,r} \over G_z} \, D_z \big ) \notag\\
  D_{p_k}G_{ij}& = E^{r,k} \Big (G_{ij,r} - \, 
{G_{,r} \over G_z}\, G_{ij,z} \Big ) \notag
 \end{align}
 
 Next  from condition G1*, we infer the existence of a mapping $X$,
 defined uniquely by
 %
 %
 \begin{equation} \label{def:X}
 {G_y \over G_z}\Big(X(y,z,q), y, z\Big) = -q
 \end{equation}
 
 \noindent for all $q\in -{G_y \over G_z}(\cdot,y,z)(\Om)$. It follows that
 %
 %
  \begin{equation} \label{def:X}
 X_q = -G_zY_p,\quad E = -G_zX^{-1}_q
 \end{equation}
 
\noindent while for fixed $y,z$ we have the following formulae for differentiation
 with respect to the $q$ variables:
%
%
 \begin{align} \label{def:D_q}
D_{q_i} &= - G_z E^{i,r}D_{x_r}  \\
 D^2_{q_\xi q_\xi} &= G^2_z E^{i,r}E^{j,s}\xi_i\xi_j D_{x_rx_s} + G_zE^{j,s}D_{x_s} 
 \big ( G_zE^{i,r} \big )\xi_i\xi_j D_{x_r} 
\notag\\
&= G^2_z E^{i,r} E^{j,s}\xi_i\xi_j  \Big \{ D_{x_rx_s}   + 
 {1 \over G_z}  \big [ E^{k,l} \big (G_zG_{rs,k} - G_{rs,z} G_{,k} \big ) D_{x_l}  \notag \\
&\qquad \qquad \qquad - \, {2 \over G_z} \, G_{s,z} D_{x_r} \big ] \Big \} \notag \\
 &= G^2_z E^{i,r} E^{j,s}\xi_i\xi_j \Big \{ D_{x_rx_s} - \big ( D_{p_l} G_{rs}\big ) D_{x_l}  - \, {2 \over G_z} \, G_{s,z} D_{x_r}  \Big \} \notag\\
&= G^2_z E^{i,r} E^{j,s}\xi_i\xi_j \Big \{ D_{x_rx_s} - \big ( D_{p_l} G_{rs} \big ) D_{x_l} \Big \} + 2E^{j,s}\xi_j G_{s,z} D_{q _\xi}, \notag
 \end{align} 
 
\noindent for any unit vector, $\xi\in\RR^n$ and $q_\xi = q.\xi$.

The proofs of Lemmas 2.1, 2.2 and 2.3 follow from formulae (2.6). First we note also that our definition of $G$-convex domain aligns with that determined by the vector field $Y$ in \cite{T1}. 
 
 \begin{lem}
 
Assume G1,G2 and G1*  hold in $\UU\times I$ with $\partial\Om\in C^2$. Then
$\Om$ is $G$-convex with respect to $y_0\in \UU^*_x, z_0\in I(\Om,y_0)$
if and only if 
%
%
\begin{equation}\label{2.7}
\left[D_i\,\gamma_j(x) - G_{ij,p_k}(x,y_0,z_0)\,\gamma_k(x)\right]
\tau_i\,\tau_j\geq 0,
\end{equation}
\noindent for all $x\in\partial\Om$,  unit outer normal
$\gamma$ and unit tangent vector $\tau$. 

\end{lem}

We prove Lemma 2.4 by applying formula (2.6) to the distance function $d = \text{dist}(\cdot,\partial\Om)$
and using the orthogonality of $ Dd = - \gamma$ and $\tau$ to cancel the last term in (2.6).

 To prove Lemmas 2.1, 2.2 and 2.3, we  let $u$ be locally $G$-convex in $\Om$ and for a
 fixed point $x_0,y_0\in \UU, z_0 = H(x_0,y_0,u(x_0))$, define the height function,
 %
 %
 \begin{equation}\label{2.8}
 h(x) = u(x) - G(x,y_0,z_0).
 \end{equation}
 \noindent Making the coordinate transformation, $x\to q = -G_y/G_z(x,y_0,z_0)$, and setting 
 $y = Tu(x) = Y(x,u(x),Du(x)), z = Z(x,u(x),Du(x))$ we compute, using formulae (2.6),
 %
 %
  \begin{align}
  D_{q_\xi}h & = -G_zE^{i, k}\xi_iD_{x_k} h,\\
 D^2_{q_\xi q_\xi}h
  &= G^2_z E^{i,k} E^{j,l} [ D_{x_kx_l}-D_{p_r} G_{kl} D_{x_r}] h \xi_i\xi_j
  + 2E^{j,s} G_{s,z} D_{q _i}h\xi_i\xi_j \notag \\
 &\ge G^2_z \big\{G_{kl}(x,y,z) -G_{kl}(x,y_0,z_0) -  D_{p_r} G_{kl}(x,y_0,z_0)\notag \\
 &\qquad \qquad [G_r(x,y,z)-G_r(x,y_0,z_0)]\big\}  E^{i, k}E^{j, l}\xi_i\xi_j 
 + 2E^{j,s} G_{s,z}\xi_jD_{q_\xi}h, \notag
   \end{align}
 
 \noindent  for any unit vector, $\xi\in\RR^n$. Setting $G_0 = G(x,y_0,z_0),
 p = G_x(x,y,z), p_0 = G_x(x,y_0,z_0)$ and using condition G4w, we  then have 
 for $h\ge 0$, that is for $u \ge G_0$
%
%
 \begin{align}
 D^2_{q_\xi q_\xi}h&\ge G^2_z [A_{kl}(x,G_0,p) -A_{kl}(x,G_0,p_0) -  D_{p_r} A_{kl}(x,G_0,p_0) \\
& \qquad \qquad (p_r-p_{0r})]E^{i, k}E^{j, l}\xi_i\xi_j + E^{j,s} G_{s,z}\xi_jD_{q_\xi}h, \notag \\
 &\ge {1\over 2} D_{p_rp_s}A_{kl}(x,u_0,p^*) D_rhD_sh(x)
                  E^{i, k}E^{j, l} \xi_i\xi_j  +  E^{j,s} G_{s,z}\xi_jD_{q_\xi}h, \notag \\  
 &\ge -K|D_{q_\xi}h|, \notag
 \end{align} 

\noindent by condition G3w and Taylor's formula, for some  $p^*$ on the straight line
segment $\ell$ joining $p$ and $p_0$ and constant K depending on  $G$,$\Om$ and $\ell$.
Setting $q_0 = q(x_0,y_0,z_0-\sigma)$, $q_t = tq + (1-t)q_0$, $x_t = X(q_t,y_0,z_0-\sigma$, $0\le t\le1$ and defining 
$h_0(t) = h(x_t))$, we can rewrite the differential inequality (2.10), 
%
%
\begin{equation}\label{2.11}
h_0^{\prime\prime} \ge -K|h_0^\prime|,
\end{equation}
\noindent which will hold whenever $h_0 \ge 0$.
For later reference we also note that $h_0^\prime(t) = D_\eta h(x_t)$ with  vector $\eta$
given by $\eta_j = E^{i,j}(q_i - q_{0i})$. 
\vspace{0.3cm}

To prove Lemma 2.1, we take $y_0 = Tu(x_0)$, which implies  that the function $G_0 = G(x,y_0,z_0)$ is a local support near $x_0$, If $h(x) < 0$, that is $u(x) < G(x,y_0,z_0)$, then we must have 
$h (x_{t_1}) > 0$ for some $t_1 \in (0,1)$, which implies that $h_0$ takes a positive maximum at 
some $t_2\in (0,1)$. Clearly this contradicts the inequality (2.11). Therefore $h(x) \ge 0$ whence $G_0$ is a global support in $\Om$ and Lemma 2.1 follows.

Note that once we have $h_0 \ge 0$, we then have that (2.11) holds everywhere and this implies a gradient estimate,
\begin{equation}\label{2.12} 
 0 \le (1-t)D_\eta h(x_t) \le Ch(x)
 \end{equation}
\noindent for some positive constant $C$, depending on $G$, $\Om$ and $Tu(\Om)$, which extends the fundamental lemma  in \cite{TW1}.

\vspace{0.3cm}

To prove Lemma 2.2, we first note that a $G$-convex function $u$ is semi-convex so that
at any singular point $x_0$, its subgradient $\partial u(x_0)$
is a convex set whose extreme points are limits of points of
differentiability. The result then follows  by showing that
the image $P_0: = G_x [(x_0,\cdot,H(x_0,\cdot,u_0)]Tu(x_0)$, where $u_0 = u(x_0)$,
is convex in $\RR^n$, (that is $Tu(x_0)$ is $G^*$-convex with respect to $(x_0,u_0)$).
Accordingly we fix two points $y_1, y_2\in Tu(x_0)$ and define corresponding
$G$-affine functions, $u_i(x) = G (x,y_i,z_i)$ for $z_i = G (x_0,y_i,u_0), i = 1,2$.
Then for any $p_0$, lying in the interior of the straight line segment joining
 $Du_1(x_0)$ and $Du_2(x_0)$,
and $y_0 = Y(x_0, u_0, p_0)$ , either $u =u_1$ or $u=u_2$ satisfies the 
 condition, $h_0^\prime(0) = D_\eta h(x_0) > 0$,
 with respect to a fixed $x$ in $Om$ or $D_\eta u_1(x_0) = D_\eta u_2 (x_0)$.
In the first cases, we obtain $h (x_{t_1}) > 0$ for some $t_1$ in $(0,1)$ close to 0 so that the proof of Lemma 2.1 is again applicable and we infer $G (x,y_0,z_0)\le \text{max}\{u_1,u_2\}(x) \le u(x)$ whenever $\eta.(Du_1-Du_2) (x_0)\ne 0$. By approximation we obtain then 
$G (x,y_0,z_0)\le u(x)$ for all $x\in\Om$ and we conclude $y_0\in Tu(x_0)$ as required.

\vspace{0.3cm}

 Lemma 2.3 follows immediately from the differential inequality (2.10), if $u\in C^2(\Om)$.
 Otherwise, we may argue along the lines of \cite{TW2,TW4}. Assuming $S_\sigma\subset\subset\Om$ and $\sigma > 0$,
 for any point $x_1\in\partial S_\sigma, u_1 = u(x_1), y_1\in Tu(x_1), z_1 = H(x_1,y_1,u_1)$,
 the inequality,
 %
 %
\begin{equation}
G(x,y_1,z_1) < G(x,y_0,z_0) + \sigma),
\end{equation}

\noindent holds for all $x\in S_\sigma$. Making the transformation
$x\to q = Q(x,y_0,z_0)$ as above, we need to show the transformed
domain $Q_0 = Q(\cdot,y_0,z_0)(\Om)$ is convex. If $Q_0$ is not convex,
there must be a straight line segment $\ell$ joining two points in $Q_0$,
containing a boundary point $q_1= q(x_1)$. Now we choose 
\begin{equation}
h(x) = G(x,y_1,z_1) - G(x,y_0,z_0)
\end{equation}
and apply the differential inequality
(2.10) along $\ell$. As before we see that h cannot take a  positive maximum on $\ell$
and hence $\ell\subset S_\sigma$. The case $\sigma = 0$ follows since
$S_0 = \cap_{\sigma>0} S_\sigma$ and Lemma 2.3 is proved.

When the function u is strictly G-convex, that is $Tu(x_0)$ is a single point for each $x\in\Om$,
then the above argument shows that $S_\sigma$ for $\sigma > 0$ is strictly $G$-convex with respect to $y_0 =Tu(x_0)$, that is the set $Q(\cdot, y_0,z_0)(S_\sigma)$ is strictly convex.

In the optimal transportation case, $G(x,y,z) = c(x,y) - z$, the c-convexity of sections is proved and used in \cite{FKM,Liu}.

To conclude this section we note that we have used condition G4w to ensure that the differential inequality (2.10) at least holds when $h \ge 0$. Otherwise it should be restricted to the set where $h = 0$.
In this case we can still prove versions of our Lemmas by strengthening other hypotheses. We will take up these results in a sequel \cite{T5}.

\vspace{1cm}

%
%

\section{Duality.}

In this section we derive the dual Monge-Amp\`ere equation for the $G$-transform $v$ of an elliptic solution u of the prescribed Jacobian equation associated with a generating function $G$ and prove the invariance of conditions G3, G3w under duality. Let $u \in C^2(\Omega)$ be elliptic for \eqref{def:MAE}   \eqref{def:GPJE} and suppose that the mapping $Tu$ defined by \eqref {def:bvp} is a diffeomorphism from to $\Omega$ to the target domain $\Omega^*$, where $G$ satisfies G1,G2,G1* and
$\Omega\times\Omega^* \subset \UU$. Then we define the local
$G$-transform of $u$ by 

%
%

 \begin{align} \label{def:vH}
  v(y) = u^*_{G,loc}(y) &= H\big (T^{-1} y, y , u\circ  T^{-1} (y) \big )\\
                                        & = Z(\cdot,u,Du)\circ (T^{-1} (y) \notag
   \end{align}

\noindent
where $H=G^*$ is the dual generating function introduced in \eqref{def:H}
and Z is defined by  \eqref{def:YZ}. When $u$ is $G$-convex in
$\Omega$, then $u^*_{G,loc}$ agrees with the G-transform defined in (1.20). 
From \eqref{def:vH} we have, for $y \in \Omega^*$,
%
%

 \begin{equation} \label{equ:Dv}
  Dv(y) = - \, {G_y \over G_z} \big ( T^{-1} y, y , v(y) \big )   
  \end{equation}

\noindent
and hence 

%
%

 \begin{equation} \label{equ:vDv}
      T^{-1} = X(\cdot , v, Dv)
  \end{equation} 
  
  \noindent
  by G1*, where $X$ is defined by \eqref{def:X}. Consequently if $u$ is an elliptic solution of the prescribed Jacobian equation \eqref{def:PJE}, \eqref  {def:example}, then $v$ is an elliptic solution of the dual prescribed Jacobian  equation,

%
%
\begin{equation} \label{def:dPJE}
\det DX (\cdot,v,Dv)  =  \psi^*(\cdot,v,Dv),\quad |\psi^*| = g/f\circ X(\cdot,v,Dv),
\end{equation}

\noindent that is

%
%
  
 \begin{equation} \label{def:dMAE} 
     {\rm det} \big [ D^2v - A^*(\cdot, v, Dv) \big ] = B^*(\cdot, v, Dv) 
 \end{equation} 
  
  \noindent
  where

%
%

 \begin{align} \label{ABstar:vDv}
     A^*(y,z,q) &= H_{yy}\left[X,y,u(X)\right] \notag \\
                       & = - \left [ \left (
       {G_y\over G_z}\right )  _ y (X,y,z) +  \left ({G_y\over G_z}\right )_z  (X,y,z) \otimes q \right ], \\
      B^* (y,z,q) &= \big|{\rm det}H_{x,y}\big|  \, {g \over f \circ X} \notag \\
                          &=
      \left ( - \, {1 \over G_z } \right ) ^n 
      \big  | {\rm det} E \, \big | 
      \, {g \over f \circ X}, \notag
  \end{align} 
 
  \vspace{0.3cm}
 
  \noindent satisfying the dual second boundary value problem:

%
%
\begin{equation} \label{def:dbvp}
T^*v(\Om^*): = X (\cdot,v,Dv)(\Om^*) = \Om.
\end{equation}

We can then formulate the dual conditions,

\vspace {0.3cm}
 
\noindent {\bf G3*}\ ({\bf G3*w}) \vspace{-0.5cm}
      $$   D_{q_iq_j} \, A_{kl}^*\; \xi _i\xi_j \, \eta_k \eta_\ell  \, > , (\ge ) \; 0, $$
 
 \vspace {0.3cm} 
  \noindent
for all $(X,y) \in \UU, z\in I(X,y) , \xi , \eta \in \mathbb{R}^{\rm n }$, $\xi \cdot \eta = 0$ . 
%
%
\begin{thm} \label{thm:condtion}
Conditions $G3, (G3w)$ and $G3^*, (G3^*w)$ are equivalent.
\end{thm}
\noindent
{\bf Proof : } It will be convenient where possible to express our formulae in terms of the matrix function $[E_{i,j}] $ defined in Condition $G2$ and the vector function $[Q_{,i}]$ given by $Q_{,i} = q_i = - \, {G_{y_i} \over G_z} $ in Condition $G1^*$. Partial derivatives of these quantities will then align with our notation in \eqref{def:Notation}.  As already indicated in Section 1, we have $Q_{i,j}= - {E_{i,j} / G_z}$. Accordingly we can express the last formula in \eqref {def:D_p}
in the form

%
%

\begin{align} \label{diff:Aij}
    D_{p_k} A_{ij} &= E^{r,k} \big (E_{ij,r} - Q_{j,r}G_{i,z} \big ) \\
    &= E^{r,k} \big (E_{ij,r} + {E_{j,r} \over G_z} \, G_{i,z} \big ) \notag \\
     &=  E^{r,k} E_{ij,r} + {G_{j,z} \over G_z} \, \delta_{i,k} . \notag 
\end{align}

\noindent
Now differentiating, with respect to $p_\ell$, we obtain

%
%

\begin{equation} \label{difftwo:Aij}
    D^2_{p_k p_\ell} A_{ij} 
    = \big ( D_{p_\ell} E^{r,k} \big )E_{ij,r}  + E^{r,k} D_{p_\ell} 
    \big (E_{ij,r}  \big )  + D_{p_\ell} \big ({G_{i,z} \over G_z} \big ) \delta _{jk} 
 \end{equation} 

\noindent
so that using \eqref{def:D_p} again we have

%
%

\begin{align} \label{diffpk:Aij}
    A^{k\ell}_{ij} &= E^{r,k} E^{s,\ell} \big \{ 
     E_{ij,rs} + E_{ij,r,z}Q_{,s} - E^{r',s'} E_{ij,r'} 
     \big ( E_{s',rs}+E_{s',r,z}Q_{,s} \big )\big \}  \\
       &\hspace{70mm} +D_{p_\ell} \big ({G_{i,z} \over G_z}\big ) \delta _{jk} \notag
    \end{align}

\noindent
To get the corresponding formula for $A^*$, we first write, 

 $$ Q_{r,k\ell} = - {1 \over G_z} \big ( E_{r,k\ell} - {E_{r,k} \over G_z} G_{, \ell,z} \big )$$
\noindent
and
\noindent
 $$ Q_{r,k,z} = - {1 \over G_z} \big ( E_{r,k,z} - {E_{r,k} \over G_z} G_{zz} \big ).$$
 
 \noindent
 Then we have

%
%

\begin{align} \label{diffpk:AijG}
    D_{q_i} A^*_{k\ell}
    &= 
    E^{i,r} \big (
     E_{r,k\ell} +  E_{r,k,z}Q_{,\ell}\big )
    - {1 \over G_z}   \big ( G_{,\ell,z}
    +G_{zz}Q_{,\ell} \big ) \delta_{ik}  
     \\
       &\hspace{70mm} +Q_{,k,z}\delta_{i\ell}.\notag
    \end{align}

\noindent
Consequently we obtain the dual formula,

%
%

\begin{align} \label{difftwopk:Aij}
   D^2_{q_i q_j} A^*_{k\ell} 
    &= -G_z
    E^{i,r}     E^{j,s}\Big \{ 
     E_{rs,k\ell} +  E_{ij,k,z}Q_{,\ell}
    - E^{r',s'} E_{rs,r'} \big(E_{s',k\ell}+E_{s',k,z}Q_{,\ell} \big) \Big \}
     \\
       &\hspace{6mm} + E^{i,r}E_{r,k,z}\delta_{j\ell}- D_{q_j} 
       \Big [ {1\over G_z} \big (G_{,\ell,z} +G_{zz}Q_{,\ell}\big ) \Big]
       \delta_{ik} + D_{q_j} \big (Q_{,k,z} \big )\delta_{i\ell} .\notag
    \end{align} 
  \vspace{.2cm} 
   
   \noindent
   Formulae \eqref{diffpk:Aij} and \eqref{difftwopk:Aij} will hold for all $(x,y,z) \in \Gamma$
   and imply the equivalence of conditions G3(G3w)  and  G3*(G3*w),  using Condition $G2$, (that is det$E \ne 0$), and the orthogonality of $\xi$ and $\eta$.
\vspace{1cm}
%
%

%
%

\section{Existence and regularity.}

We begin with the definition of generalized solution. Let $\Omega$ and $\Omega^*$ be bounded domains in $ \RR^n$, with $\overline\Omega\times\overline\Omega^*\subset\UU$, and let 
$u\in C^0(\overline\Omega)$ be $G$-convex in  $\Omega$, with conditions G1,G2, G1* satisfied. Following \cite{MTW}, we introduce a measure $\mu= \mu_g[u]$ on $\Omega$, for $g \ge 0 \in L^1(\Omega)$, such that for any Borel set $E \subset \Omega$, 

%
%
\begin{equation} \label{def:measmu}
\mu (E) = \int _{Tu(E) } g
\end{equation}

To prove that $\mu$ is a Radon measure, we can extend the argument in \cite{MTW}, using the $G$-transform  $v$, defined in \eqref {def:funcv}, namely

$$v(y) = u^*_G(y) = \text{sup}_{\Om} H(\cdot,y,u).$$

\vspace {0.2cm}

\noindent  Note that v will be $G^*$- convex in $\Omega^*$.
 Defining the dual $G^*-$transform on $\Omega^*$, $v^*$, by

%
%
\begin{equation} \label{def:vstar}
  v^*(x) = \sup_{y\in \Omega^*} G(x,y,v(y)), 
\end{equation}

\noindent
we obtain $v^*= u $ in $\Omega$ and furthermore $y \in Tu(x)$ if and only if $x\in T^*v(y)$, where $T^*$ denotes the $G^*$-normal mapping on $\Omega^*$. Since the functions $u$ and $v$ are semi-convex and hence twice differentiable almost everywhere in $\Omega$ and $\Omega^*$, respectively, we then infer that the subset of points $y \in \Omega^*$, such that $y \in Tu(x_1)\cap Tu(x_2)$ for some $x_1 \ne x_2$, $\in \Omega$, has measure zero and from this it follows that $\mu$ is countably additive. Also extending the argument in \cite{MTW}, we have that $\mu$ is weakly continuous with respect to local uniform convergence, namely if $\{u_m\}$ is a sequence of $G$-convex functions in $\Omega$, converging to $u$, then the sequence of measures $\{\mu _g[u_m]\}$ converges to $\mu_g[u]$ weakly. These arguments also parallel the special case of the Monge-Amp\`ere measure of Aleksandrov, as presented for example in the book \cite{Gu}.

A $G$-convex function $u$ on $\Omega$ is now defined to be a generalized solution of the second boundary value problem \eqref{def:bvp} for equation \eqref{def:PJE}, \eqref{def:separable}, under the mass balance condition \eqref{def:massbalance}, if

%
%
\begin{equation} \label{equ:mug}
  \mu_g[u] = \nu_f
\end{equation}

\vspace{0.2cm}

\noindent where $\nu_f = f {\rm d}x$ and $g$ is extended to vanish outside $\Omega^*$. More generally we can replace $\nu_f$ by any finite Borel measure $\nu$ on $\Omega$. We remark that this notion corresponds to that of generalized solution of Type A in \cite{KW}, where a generalized solution of Type B corresponds to the dual notion, 
%
%
\begin{equation} \label{eqn:muf}
\mu^*[u] (E^*) := \int _{T^{-1}(E^*)} f = \nu^*(E)
\end{equation}

\noindent for any Borel set $E^* \subset \Omega^*$, that is $\mu_f[v] = \nu^*$, where $\nu^*=g\,{\rm d}y$.  As in the optimal transportation case in \cite{MTW}, the two notions are equivalent provided the measures $\nu$ and $\nu^*$ have densities $f \in L^1(\Omega)$, $g\in L^1(\Omega^*)$ respectively. Therefore we have:

 \begin{lem} 
Let u be a generalized solution of the second boundary value problem \eqref{def:bvp} for equation \eqref{def:PJE}, \eqref{def:separable}, under the mass balance condition \eqref{def:massbalance},
and let v be the G-transform of u. Then v is a generalized solution of the second boundary value problem  \eqref{def:dPJE}, \eqref{def:dbvp}.
\end{lem}

Next we formulate an existence theorem under condition G5 which extends the optimal transportation case.

\begin{thm}
Let $\Omega$ and $\Omega^*$ be bounded domains in $ \RR^n$, with $\overline\Omega\times\overline\Omega^*\subset\UU$ and let G be a generating function satisfying G1,G2,G1* and G5.
Suppose that f and g are positive densities in $L^1(\Om)$ and $L^1(\Om^*)$ satisfying the mass balance condition \eqref{def:massbalance}.Then for any $x_0\in\Om$ and $u_0 > m_0 +K_1$,
where $K_1 = K_0{\rm dist}(x_0,\partial\Om)$, there exists a generalized solution of \eqref{def:MAE},\eqref{def:bvp} satisfying $u(x_0) = u_0$.
\end{thm}

The proof of Theorem 4.2 follows the approach in \cite{CO,KO} using approximations by piecewise $G$-affine functions to solve the dual problem \eqref{eqn:muf}, where $\nu^*$ is approximated by a linear
combination of Dirac delta measures. Condition G5 implies that a generalized solution $u$ whose graph passes through $(x_0,u_0)$ satisfies $ u \ge m_1$ for some $m_1 > m_0$, together with an a priori bound:
%
%
\begin{equation}
 |u|_{0,1} \le K,
 \end{equation}
 \noindent with constant $K$ depending on $K_0, u_0$ and diam $\Om$. The solvability of the dual problem, with an arbitrary finite measure $\nu^*$ then follows from the weak continuity of $\mu^*$.
Theorem 4.2 then follows from Lemma 4.1.

Taking account of our remark following the formulation of G5 in Section 1 and choosing $\mu(G) = 
\text{log}(G)$ with $m_0 = 0$, we see that if we replace the gradient bound in G5 by 
%
%
\begin{equation}
|G_x(x,y,z)| \le K_0 G (x,y,z),
\end{equation}

\noindent we conclude that there exists a generalized solution $u$ such that $u(x_0) = u_0$ for any given $u_0 > 0$. We remark also that instead of prescribing $u(x_0)$ in Theorem 4.3, 
we may fix instead $u_{min} = m_1 > m_0$.

The theory developed in the preceding sections can now be applied to prove interior regularity results for generalized solutions. A crucial lemma for this purpose is Lemma 2.2 on the characterization of the mapping $Tu$ by the subdifferential $\partial u$, which shows that the concept of generalized solution is local, that is it is the same for subdomains. This property was originally overlooked in \cite{MTW} but rectified subsequently in \cite{TW4}, with a more direct approach found in \cite{KM}. In this paper we will just prove interior smoothness under condition G3 following the approach in \cite{MTW}. Further extensions of optimal transportation regularity depend also on Lemma 2.3 on the $G$-convexity of sections but we will not pursue these here. First we state some more lemmas which extend the corresponding results in \cite{MTW}.

\begin{lem}
Let u be a generalized solution of \eqref{def:MAE}, \eqref{def:GPJE},\eqref{def:bvp}. 
Suppose $f > 0$ in $\Om$ and $\Om^*$ is $G$-convex with respect to $u$. 
Then $Tu(\Om) \subset \overline\Om^*$.
\end{lem}

\begin{lem}
Let $\Om$ be a bounded domain in $\RR^n$ and $u,v \in C^0(\overline\Om)$, $G$-convex functions satisfying $u \ge v$ in $\Om$ and $u = v$ on $\partial\Om$. Then $Tu(\Om) \subset Tv(\Om)$.
\end{lem}

The proofs of Lemmas 4.3 and 4.4 are essentially identical with the optimal transportation case in \cite{MTW}. Lemma 4.3 makes critical use of the subdifferential being convex in general.
 To prove Lemma 4.4 we fix some $x_0\in \Om$, $y_0 \in Tu(x_0)$ and increase $z_0 = H(x_0,y_0, u(x_0))$ in $I$ until the function $G(\cdot, y_0, z_1)$, for $z_1 > z_0$, becomes a $G$-support for $v$ at some point $x_1\in\Om$. But then, using Lemma 2.2 if $v$ is not differentiable at $x_1$, we must have $y_0\in Tu(x_1)$ since $v$ is $G$-convex.

\vspace {0.2cm}

Now we formulate the interior regularity result, which extends the main result in \cite{MTW}.

\begin{thm}
Let $u$ be a generalized solution of \eqref{def:bvp} with positive densities $f \in C^{1,1}(\Om), g \in C^{1,1}(\Om^*)$ with $f,1/f \in L^\infty(\Om), g,1/g \in L^\infty(\Om^*)$ and with generating function G satisfying conditions, G1,G2,G1*, G3 and G4w. 
Suppose that $\Om^*$ is $G^*$-convex with respect to $u$. Then $u\in C^3(\Om)$ and is an elliptic solution of \eqref{def:MAE}, \eqref{def:GPJE}. Furthermore if $\Om$ is $G$-convex with respect to $v = u^*_G$, then $Tu$ is also a diffeomorphism from $\Om$ to $\Om$*, with $v$ an elliptic solution of the dual equation \eqref{def:dMAE}.
\end{thm}

To prove Theorem 4.5 by the method in \cite{MTW}, we need to solve approximating Dirichlet problems,
%
%

\begin{equation} \label{def:DP}
\begin{align}
\det [D^2w - A(\cdot,w,Dw) ]& =  B(\cdot,w,Dw) \quad\mbox{in }B_r,\\
                                              w& = u_m  \qquad\mbox{ on }\partial B_r, \notag
\end{align}
\end{equation}

\noindent where $B_r$ is a small ball in a fixed subdomain $\Om^\prime \subset\subset\Om$ and $\{u_m\}$ is a sequence of smooth functions converging uniformly in $\Om^\prime$ to $u$. 
Specifically, we set

$$ u_m=\tilde u_{h_m}- c_0|x|^2,$$

\noindent where $\tilde u = u+ c_0|x|^2$ is convex in $\Om^\prime$ and $\tilde u_{h_m}$ denotes the mollification of $\tilde u$ with $h_m \to 0$. It follows also from the local Lipschitz continuity and semi-convexity of the G-convex function $u$, then that for $m$ sufficiently large, the constant
 $c_o$ can be chosen so that

$$ (|u_m| + |Du_m|) \le c_o, \ \ D^2u_m \ge -c_0I $$

\noindent in $\Om^\prime$ and moreover $u_m$ is admissible in the sense that 
%
%
\begin{equation} \label{constraint}
(x,u_m(x), Du_m(x)) \in \Gamma^\prime(\Om)
\end{equation}

\noindent for all $x\in \Om^\prime$, where in accordance with our notation in Section 1,

 $$\Gamma^\prime(\Om) = \big \{(x,u,p) \mid x\in\Om, u = G (x,y,z), p = G_x(x,y,z), 
\text {for some} \ y\in\UU^*_x, z\in I(x,y)\}.$$

 The essential  difference here with the corresponding equations in
\cite{MTW} is the dependence of the matrix function $A$ and scalar function $B$ on $u$. The necessary existence and uniqueness result is expressed in the following lemma.

\begin{lem}
Let G satisfy conditions G1, G2, and G3w in $\UU$, $\Om\subset \UU^{(1)}$ and suppose 
$B > 0, \in C^{1,1}(\Gamma^\prime)$. Let $B_r$ be a small ball in
 $\Om^\prime\subset\subset\Om$ and suppose $u_m\in C^4(\Om^\prime)$ satisfies the constraint  \eqref {constraint}.
Then for sufficiently small $r$, depending on $G , B$ and $c_o$,
there exists a unique elliptic solution $w\in C^3(\bar B_r)$
of the Dirichlet problem  \eqref{def:DP}, also satisfying \eqref{constraint} in $B_r$
 together with an a priori $C^1$ bound,
%
%
\begin{equation}
 |w| + |Dw| \le C,
 \end{equation}
\noindent with constant  $C$ also depending only on $G, B$ and $c_0$.
\end{lem}

To prove Lemma 4.6 by the method of continuity, (as in \cite{GT}), we need suitable a priori estimates for solutions and their derivatives up to second order. For simplicity we can
assume $B_r = B_r(0)$ is centred at the origin.  First we note that, even with
 the dependence of $A$ and $B$ on the solution $w$, the functions 
 
     $$  v=v_m=u_m+k(|x|^2-r^2)$$
     
\noindent will still be strict elliptic upper barriers for \eqref{def:DP} for sufficiently large $k$ and small $r$, depending on $A,B$ and $c_0$, as well as satisfy \eqref{constraint} in $B_r$. 
 Moreover for any constant $c_1$, we can obtain by such choice the strong differential inequality,
 
 $$\det [D^2v - A(\cdot,v,Dv) ] >  B(\cdot,v,Dv) + c_1\quad\mbox {in }B_r.$$
 \vspace{0.2cm}

\noindent Now let w be an elliptic solution of \eqref{def:DP} satisfying \eqref{constraint} in $B_r$ and fix $k=k_0$ and  $r=r_0$ as above. We claim then that for $k = k_0$ and sufficiently smaller $r < r_0$, we have $w \ge v$ in $B_r$. 
To show this we suppose $M = \text{max}(v-w) > 0$ is taken on at 
 $x_0 \in \Om$. Since $Dv(x_0) = Dw(x_0)$ and $D^2v(x_0) \le D^2w(x_0)$ we obtain, by increasing $k$ appropriately,
 
 $$ w-v \le CMr^2$$
 
\noindent  for further constant $C$ depending on $A,B, k_0, r_0$ and $c_0$. Choosing $r$ sufficiently small,we infer  $w \ge v$, as claimed. We obtain thus an estimate for $w$ from below.
Furthermore we can  also obtain, since $v$ is convex for $k\ge c_o$,
%
%
\begin{equation}
 w > v, \ \ Dw (\partial B_r) \subset Dv (\bar B_r)
 \end {equation}
 
 \vspace{0.2cm}
 
 \noindent which provides an estimate for $Dw$ on $\partial B_r$.
 Next  since the Jacobian determinant, det$DTw\ne 0$ in $B_r$, we see that $|Tw|^2$ takes its maximum on $\partial B_r$ so that we obtain an estimate in $B_r$ for $Tw$. 
 To complete the gradient bound (4.9), we need an estimate for $w$ from above. Here we use the ellipticity to obtain 

$$ \triangle w \ge \text{traceA} (\cdot.w,Dw)$$
\vspace{.2cm}

\noindent so that with a similar argument we can show 

$$w\le \underline v = \text{sup}u_m - k(|x|^2 - r^2)$$

\vspace{.2cm}

\noindent in $B_r$ for sufficiently large $k$ and small $r$. From the estimates for $w$ and $Tw$, we then have an estimate for $Dw$ in $B_r$ and so (4.9) is proved.

We remark that if G1* and G4 are also satisfied then both $w$ and $v$ will be $G$-convex in
 $B_r$ by Lemma 2.1 and we can use Lemma 4.4 to obtain $Tw (B_r) \subset Tv (B_r)$ and hence alternatively estimate $Tw$ in $B_r$.
 
Global second derivative estimates for solutions of \eqref{def:DP}  follow as in the case
 where there is no dependence on $u$
in $A$ and $B$, since  functions $\varphi$ given by
%
%
\begin{equation}
\varphi (x) =  k(|x|^2 - r^2)
\end {equation}

\noindent provide barriers for the linearized operators in $B_r$  for sufficiently large $k$ and small $r$, that is
%
%

\begin{equation}
[D_{ij}\varphi -
D_{p_k}A_{ij}(\cdot,w,Dw)D_k\varphi]\xi_i\xi_j \ \ge \  |\xi|^2
\end{equation}
\vspace{0.2cm}

\noindent in $B_r$, for all $\xi\in \RR^n$. The reader is referred to the papers \cite {JTY,T1,TW3} for more details. 

The uniqueness of solutions of \eqref{def:DP} can be shown by the same comparison argument as above or by using the linearized equation and (4.12). Once second derivative bounds and uniqueness are established, the standard method of continuity in \cite{GT} is applicable  with the smallness of the radius r also used to imply the invertibility of the linearized operators in the associated deformation. 

In order to proceed from Lemma 4.6 we  send $m$ to $\infty$, that is for $u_m$
to approach $u$. At this stage we need to replace G3w by the strict regularity condition G3 to obtain an interior second derivative estimate for solutions, namely for any $r^\prime < 1$
%
%

\begin{equation}
|D^2w| \ \le \ C
\end{equation}

\noindent in $B_{r^\prime}$ where $C$ depends on $n,A,B, \text{sup}_{B_r}(|w| + |Dw|)$
 and $r - r^\prime$. 
In the optimal transportation case this was the key estimate in [MTW] and as remarked there it 
also embraces general equations of the form \eqref{def:MAE}; (see \cite{TW2} for a direct proof).

Using also the gradient estimate (4.9) we  then conclude that the solutions $w_m$ of
 \eqref{def:DP} converge to a unique solution $w\in C^3(B_r)\cap C^{0,1}(\bar B_r)$
 of the limiting problem $w = u$ on $ \partial B_r$. To complete the proof of Theorem 4.5, we need to show that $w = u$.  This is more delicate than in the optimal transportation case as we cannot localise around maximum or minimum points  but we overcome this obstacle by  using both Lemmas 2.1 and 2.2. Let us first suppose that $w - u$ takes a positive maximum $M$ at some point $x_0$ in $B_r$. Since $u$ and $w$ lie in $C^{0,1}(\bar B_r)$ there exists a radius $r^\prime < r$, depending on $M$ and the Lipschitz constants of $u$ and $w$, such that $w-u < M/2$ on 
 $\partial B_{r^\prime}$. Since $w \in C^2(\bar B_{r^\prime}$, for sufficiently small
 $\epsilon > 0$, we can perturb $w$ to get a  function $w_\epsilon \in C^2(\bar B_{r^\prime})$, satisfying $|w - w_\epsilon| < \epsilon$, which is a strict upper barrier for (4.7), 
 that is
%
%
 \begin{equation}
 \det [D^2w_\epsilon - A(\cdot,w_\epsilon,Dw_\epsilon) ] >  
 B(\cdot,w_\epsilon,Dw_\epsilon) \end{equation}
 \vspace{0.2cm}
 
  \noindent in $B_{r^\prime}$. Furthermore $w_\epsilon$ is elliptic for (4.7) and 
  satisfies (4.8) in $B_{r^\prime}$. Hence for small enough $\epsilon < M/4$ the subdomain
  $\Om_\epsilon = \{ w_\epsilon > u\} \subset B_{r^\prime}$ and moreover by Lemma 2.1,
   $w_\epsilon$ is $G$-convex in $\Om_\epsilon$ so that  from (4.14), 
  %
  %
   \begin{align}
    \mu_g[w_\epsilon] (\Om_\epsilon) & > \int_{\Om_\epsilon} f \\
    \vspace{.2cm}
                                                                  & =  \mu_g[u] (\Om_\epsilon) \notag
     \end{align}                                                                                                                             
   
    \noindent  by Lemma 2.2.                                                               
 Now using the monotonicity, Lemma 4.3, we reach a contradiction and hence $w \le u$ in $B_r$.
 Similarly we show $w\ge u$ in $B_r$, whence $w = u$ in $B_r$ and consequently
  $u \in C^3(\Om)$ as asserted.  

We point out however that we have implicitly used a stronger condition on the target density $g$ in the proof of Theorem 4.5, namely $g \in C^{1,1}(\overline\Om^*)$. To use only the local smoothness of $g$ we need to prove in advance that $u\in C^1(\Om)$, (to keep $T(B_r)$ away from $\partial \Om^*$), or equivalently by duality that $u$ is strictly $G$-convex.This may be accomplished under the hypotheses  $f, 1/f \in L^\infty(\Om), g,1/g \in L^\infty(\Om^*)$ as in \cite{TW3} or \cite{Loe} and will be treated in \cite{T5} in conjunction with the relaxation of condition G4w. Moreover, as in \cite{Loe}, we
also infer $u\in C^{1,\alpha}(\Om)$ for some $\alpha > 0$.

Under our smoothness conditions on f and g in Theorem 4.5, we actually have
the solution $u \in C^{3,\alpha}$ for all $\alpha < 1$, by the Schauder theory \cite{GT}.
We also obtain that if $f \in C^\infty(\Om), g \in C^\infty(\Om^*)$, then 
$u\in C^\infty(\Om)$. 
As in the case with optimal transportation, regularity results may be refined when
 perturbation arguments using Lemma 2.3 are employed instead of Lemma 2.2. In particular,
when we replace $C^{1,1}$ by $C^{0,\alpha}$ in Theorem 4.5, we obtain from \cite{LTW} that $u\in C^2(\Omega)$. We will not pursue this approach in this paper. It would be interesting though to have extensions to G3w along the lines of \cite{FKM,Vet}, that is to show solutions are strictly
G-convex and continuously differentiable under G3w, together with appropriate domain convexity conditions. We remark also that using the interior second derivative estimates in \cite{LT1}, it follows from our  proof of Theorem 4.5 that condition G3 in the hypothesis can be relaxed to G3w if the solution $u$ is assumed strictly $G$-convex in $\Om$.

We also remark here that if condition G3w is violated then Lemma 2.2 is no longer true in general and that generalized solutions are not necessarily continuously differentiable. This follows in the same way as the optimal transportation case in \cite{Loe, MTW}. Similarly from \cite{MTW}, the $G^*$- convexity of the target
$\Om^*$ is also necessary for regularity. An interesting consequence is  that if the generating function satisfies G3, then generalized solutions of the complementary problem associated with condition (1.16)  cannot be continuously differentiable in general.

Combining Theorems 4.2 and 4.5 we get an existence result for locally smooth solutions.

\begin{cor}
Let $\Omega$ and $\Omega^*$ be bounded domains in $ \RR^n$, with $\overline\Omega\times\overline\Omega^*\subset\UU$ and let G be a generating function satisfying G1,G2,G1*,G3,G4w and G5.
 Suppose that $f \in C^{1,1}(\Om), g \in C^{1,1}(\Om^*)$ with $f,1/f \in L^\infty(\Om), g,1/g \in L^\infty(\Om^*)$ and that $f$ and $g$ satisfy the mass balance condition  \eqref{def:massbalance}.Then for any $x_0\in\Om$ and $u_0 > m_0 +K_1$, there exists a $G$-convex, elliptic solution $u\in C^3(\Om)$of \eqref{def:MAE},\eqref{def:bvp} satisfying $u(x_0) = u_0$ and $Tu(\Om) = \Om^*$ a.e., provided $\Om^*$ is
 $G^*$-convex with respect to $\Om\times(m_0,\infty)$. If also $\Om$ is $G$-convex with respect to
 $\Om^*\times I$, then $Tu$ is a diffeomorphism from $\Om$ onto
 $\Om^*$.
 \end{cor}
 
 We remark that the interval $(m_0.\infty)$ in the $G^*$-convexity hypothesis may be replaced by the subinterval, $(u_0 - K_1,u_0+K_1)$, while to ensure $Tu$ is onto $\Om^*$ we need only assume that $\Om$ is $G$-convex with respect to all $y\in\Om^*, z\in I $ satisfying
 
 $$ |G(\cdot,y,z) - u_0| <K_1. $$ 

\vspace{0.3cm} 

\noindent
 Furthermore  if (4.6) holds then we can take any $u_0 > m_0 = 0$ in Corollary 4.7. In \cite{T5}, we show that condition G4w is not needed for Corollary 4.7 as the other hypotheses there will still enable the applications of Lemmas 2.1 and 2.2.
 
 Theorems 4.2,4.5 and Corollary 4.7 clearly apply to the parallel beam example (1.18) and under the hypotheses of Corollary 4.7 we would obtain the existence of a positive elliptic solution $u$ of equation (1.18) with $|Du| <1$. 
  
 We conclude this section with some remarks about the point source reflection problem treated in 
 \cite{KO,KW,LT2}. Invoking the ellipsoid of revolution in its polar form  \cite{KO,KW} and projecting
 onto the unit ball in $\RR^n$,  we may take $\UU = B_1(0)\times \RR^n, I = (0,\infty)$. In particular, if the target hypersurface lies in a graph $y_{n+1} = T(y)$ in $\RR^{n+1}$, we may model the reflection process by  a generating function:
 %
 %
 \begin{equation} \label{pointsource}
 G(x,y,z) = {1\over z}\big\{\sqrt{z + |y|^2 + T{^2}} - x.y - \sqrt{1-|x|^2}\ T\big\}.
 \end{equation}
 
 \noindent Differentiating we obtain,
 %
 %

 \begin{align}
  G_x  
  & = {1\over z} (-y + {x\over\sqrt{1-|x|^2}}\ T), \\
 G_{xx} 
 & = {1\over z(1-|x|^2)^{3/2}} (I + x\otimes x) T. \notag
  \end{align}

 \vspace{0.3cm} 
 
  \noindent Let us restrict here to a special case where the target lies in a hyperplane not above the source at the origin, namely $y_{n=1} = \tau \le 0$ . Writing $\bar u = u - p.x, p = G_x$,  it follows that  $\bar u > |p|$ and moreover we have the formulae:
  %
  %
  \begin{align}
  Z  & = (1 - {2\tau u\over\sqrt{1-|x|^2}}) / ({\bar u}^2 - |p|^2), \\
  Y & = - pZ + {\tau x\over \sqrt{1-|x|^2}}, \notag \\
  A & = {\tau ({\bar u}^2 - |p|^2)\over(1-|x|^2)(\sqrt{1-|x|^2} - 2\tau u)}(I + x\otimes x). \notag
  \end{align}
  
  It is now easy to check that $G$ satisfies G1,G2,G1*,G2,G3,G4 with $I(x,y) = J(x,y) = (0,\infty)$ if $\tau < 0$, and G1,G2,G1*,G3w,G4w for $\tau = 0$, in which case $A = 0$. Furthermore if $\bar\Om \in B_1(0)$, then $G$ also satisfies (4.6). The resultant existence and regularity results correspond to special cases of  \cite{KW} except we are also able to solve the boundary value problem (1.7) in the classical sense when $\Om$ is also G-convex. However the local regularity results in \cite {KW} are more general in that
 condition G3 need not be satisfied everywhere for more general targets and in this case a stronger version of Lemma 2.2 enables local regularity to be proved  in the sub-domains where G3 is satisfied. When $\tau = 0$, $G$-convexity coincides with convexity and local regularity under $G^*$-convexity  of the target domain $\Om^*$ follows from \cite{C1}, Lemma 3.
 
 Finally we note that the existence of globally smooth solutions of the point source and parallel beam near field reflection problems is established in \cite{LT2} and that solutions of these and related problems can also be represented as potentials of a nonlinear Kantorovich problem, \cite{GO,Liu2,Liu3}.

\vspace{1cm}
%
%
%
%


\frenchspacing
\begin{thebibliography}{99}


\bibitem{C1} L. Caffarelli,  The regularity of mappings with a convex potential,
\emph{J.\ Amer.\ Math.\ Soc.\ }\textbf{5} (1992), 99--104.

\bibitem{C2} L. Caffarelli,  Boundary regularity of maps with convex potentials
II. \emph{Ann.\ of Math.\ }\textbf{144}  (1996), 453--496.

\bibitem{CO} L. Caffarelli and V.I. Oliker,  Weak solutions of one inverse problem in geometric optics. \emph{ J. Math. Sci. (N. Y.)} \textbf{154} (2008),  39--49.

\bibitem{D} Ph. Delano\"e, Classical solvability in
dimension two of the second boundary value problem associated with
the Monge-Amp\`ere operator, \emph{Ann.\ Inst.\ Henri Poincar\'e,
Analyse Non Lin\'eaire} \textbf{8} (1991), 443--457.

\bibitem {FKM} Figalli, A., Kim, Y.-M., McCann, R.,
             H\"older continuity and injectivity of optimal maps, 
               arXiv:1107.1014. 
               

\bibitem{GT} D. Gilbarg, and N.S. Trudinger,
\emph{Elliptic partial differential equations of second order}, Springer, 1983.


\bibitem{GO} T. Graf and V. Oliker, An optimal mass transport approach to the near-field reflector problem in optimal design, 
\emph{Inverse Problems} \textbf{28} (2012), 1-15.


\bibitem{GH} C.E. Gutierrez and Q. Huang, The near field refractor, preprint (2011).

\bibitem{Gu} C.E. Gutierrez, \emph {The Monge -Amp\`ere equation}, Birkh\" auser, 2001.

\bibitem{JTY} F. Jiang, N. S. Trudinger and X.-P. Yang,
On the Dirichlet Problem for Monge-Amp\`ere type equations,
\emph {Calc. Var. Partial Differ. Equ.}, to appear, online DOI 10.1007/s00526-013-0619-3.

    \bibitem{KM}
 Y.-H. Kim and R.J. McCann,
  Continuity, curvature and the general covariance of
    optimal transportation,
    \emph{J. Eur. Math. Soc.} \textbf{12} (2010), 1009--1040.

   \bibitem{KO}
       S.A. Kochengin and V.I. Oliker,
    Determination of reflector surfaces from near-field
    scattering data,
    \emph{Inverse Problems} \textbf{13} (1997), 363-373.

\bibitem{KW}
     A. Karakhanyan and X.-J. Wang,
    On the reflector shape design,
     \emph{J. Diff. Geom.}  \textbf{84} (2010), 561--610.

\bibitem{LTU}  P.-L. Lions, N.S.Trudinger and J. Urbas,
Neumann problem for equations of Monge-Amp\`ere
 type,
  \emph{Comm.\ Pure\  Appl.\ Math.\ }\textbf{39} (1986), 539--563,
  \emph{Proc.\ Centre\ Math.\ Anal.\ }\textbf{10} (1985), 135--140.

\bibitem{Loe}
     G. Loeper,
    On the regularity of solutions of optimal transportation problems,
     \emph{Acta Math.} \textbf{202} (2009), 241--283.
    
  \bibitem{Liu}
     J.-K. Liu,
     H\"older continuity of optimal mappings in optimal transportation,
     \emph{Calc.\ Var.\ PDE} \textbf{34} (2009), 435--451.
    
  \bibitem{Liu2}
      J.-K. Liu,
     Light reflection is nonlinear optimization,
     preprint (2012).   
    
  \bibitem{Liu3}
    J.-K. Liu,
     A class of nonlinear optimization problems with potentials,
     preprint (2012).
 
 \bibitem{LT1}
      J.-K. Liu and N.S Trudinger
      On Pogorelov estimates for Monge-Amp�re type equations,
      \emph{Discrete Contin. Dyn. Syst}. Ser. A. \textbf{28}, (2010), 1121--1135.

 \bibitem{LT2}
      J.-K. Liu and N.S Trudinger,
    On classical solutions of near field reflection problems,
    preprint (2013).
    
 \bibitem{LTW} J.-K. Liu, N.S Trudinger and X.J. Wang, 
       Interior $C^{2,\alpha}$ regularity for potential functions
       in optimal transportation.
       \emph{Comm. Partial Differential Equations} \textbf{35} (2010), 165--184.

\bibitem{MTW} X.-N. Ma, N.S.Trudinger and X.-J.Wang,
              Regularity of potential functions of the optimal
              transportation problem,
             \emph{Arch.\ Rat.\ Mech.\ Anal.\ } \textbf{177} (2005), 151-183

\bibitem{RR} S.T. Rachev and L.  Ruschendorff, \emph{Mass transportation problems},
Springer, 1998.

\bibitem{T1}  N.S.Trudinger,
Recent  developments in elliptic partial differential equations
of  Monge-Amp\`ere type, \emph {Proc.\ Int.\ Cong.\ Math.},Madrid,
\textbf{3} (2006), 291-302.

\bibitem{T2}
     N.S.Trudinger, 
   On the prescribed Jacobian equation,
    Gakuto Intl. Series, Math. Sci. Appl. \textbf{20} (2008),
    Proc. Intl. Conf. for the 25th Anniversary of Viscosity Solutions, 243--255
 
  \bibitem{T3}
    N.S.Trudinger, 
   On the local theory of prescribed Jacobian equations,
    Conference in honor of 60th birthday of C.-S. Lin, Taipei, 2011,
                                (https://maths.anu.edu.au/~neilt/RecentPapers.html).
 
  \bibitem{T4}
    N.S.Trudinger, 
     On generated prescribed Jacobian equations,
     \emph{Oberwolfach Reports} \textbf{38} (2011), 32-36.
     
   \bibitem{T5}
    N.S.Trudinger,
   The local theory of prescribed Jacobian equations revisited, (in preparation).
    

\bibitem{TW1}
     N.S.Trudinger and X.-J. Wang,
     The Monge-Amp\`ere equation and its geometric applications,
     in \emph{Handbook of geometric analysis}, International Press (2008)
     467-524.



\bibitem{TW2}
    N.S.Trudinger and X.-J. Wang,
      On convexity notions in  optimal transportation,
     preprint  (2008).

\bibitem{TW3}
    N.S.Trudinger and X.-J. Wang ,
     On the second boundary value problem for Monge-Amp\`{e}re type equations and optimal transportation,
   \emph{Ann. Scuola Norm. Sup. Pisa Cl. Sci.},Series 5, \textbf{8}, (2009), 143--174.

\bibitem{TW4}   N.S.Trudinger and X.-J.Wang,
              On strict convexity and continuous differentiability of potential
              functions in optimal transportation,
               \emph{Arch.\ Rat.\ Mech.\. Anal.\ }, \textbf{192} (2009), 403--418.

    
  \bibitem{U} J. Urbas,  On the second boundary value problem for equations of
Monge-Amp\`ere type, \emph{J.\ Reine Angew.\ Math.\ }\textbf{487} (1997),
115--124.

           
 \bibitem{Vet} J.Vetois,  Continuity and injectivity of optimal maps, preprint (2011).


\bibitem{V1} C.Villani,  \emph{Topics in optimal transportation},
Amer.\ Math.\ Soc., 2003.

\bibitem {V2} C.Villani,
             \emph{Optimal transportation, old and new},
             Springer, 2008
.


\end{thebibliography}
\end{document}